\documentclass[10pt]{article}

\usepackage{amsfonts}
\usepackage{amsmath,amssymb,amsthm}
\usepackage{url}\usepackage{amsmath,amssymb,amsthm}
\usepackage{url}
\usepackage{epsf,epsfig,amsfonts,a4wide}
\usepackage{epsf,epsfig,amsfonts}

\usepackage{amsfonts}
\usepackage{amsmath,amssymb,amsthm}
\usepackage{url}\usepackage{amsmath,amssymb,amsthm}
\usepackage{url}\textwidth=1.06\textwidth

\newtheorem{rem}{Remark}
\newtheorem*{Tha}{Theorem A}
\newtheorem*{Thb}{Theorem B}

\newtheorem*{Lemma}{Lemma}

\theoremstyle{remark}

\newcommand{\const}{\mbox{\rm const}}

\newcommand{\weg}[1]{}

\begin{document}
\title{Appendix: Dini theorem for pseudo-Riemannian metrics }
\author{
Alexei V. Bolsinov\thanks{	
Department of Mathematical Sciences,  
Loughborough University,  LE11 3TU 
UK, A.Bolsinov@lboro.ac.uk}, 
Vladimir S. Matveev\thanks{Institute of Mathematics, FSU Jena, 07737 Jena Germany,  matveev@minet.uni-jena.de}, 
Giuseppe Pucacco\thanks{Dipartimento di Fisica, Universit\`a di Roma ``Tor Vergata", 00133 Rome Italy, pucacco@roma2.infn.it}}
\date{}

\maketitle

\section{Introduction} 

Consider a  Riemannian or a  pseudo-Riemannian metric $g=(g_{ij})$ on  a surface $M^2$.  We say that a  metric $\bar g$ on the same surface is \emph{projectively equivalent} to $g$,
 if every geodesic of $\bar g$ is a reparametrized geodesic of $g$.   
In 1865 Beltrami  \cite{Beltrami} asked\footnote{ Italian original from \cite{Beltrami}: 
La seconda $\dots$  generalizzazione $\dots$ del nostro problema,   vale a dire:   riportare i punti di una superficie sopra un'altra superficie in modo  che alle linee geodetiche della prima corrispondano linee geodetiche della seconda.}  to describe all pairs of  projectively equivalent Riemannian  metrics on  surfaces. 
From the context it is clear that he considered this problem locally, in a neighbourhood of almost every  point.  

Theorem A below, which is the main result of this note,  gives  an answer to the following generalization of  the question of Beltrami: we allow the metrics $g$ and $\bar g$ to be pseudo-Riemannian.

\begin{Tha} \label{thm3}
Let $g$, $\bar g$ be projectively equivalent  metrics on $M^2$, and $\bar g\ne \textrm{const} \cdot g$ for every $\textrm{const}\in \mathbb{R}$. Then, in the neighbourhood of almost every point there exist coordinates $(x,y)$  such that the metrics are  as in the following table. \\
 \begin{tabular}{|c||c|c|c|}\hline &  \textrm{Liouville case} & \textrm{Complex-Liouville case} & \textrm{Jordan-block case}\\ \hline \hline
$g$ & $(X(x)-Y(y))(dx^2 \pm  dy^2)$ &  $2\Im(h)dxdy$ & $\left( 1+{x} Y'(y)\right)dxdy $
\\  \hline  $ \bar g$  &$ \left( \frac{1}{Y(y)}-\frac{1}{X(x)}\right) \left( \frac{dx^2}{X(x)} \pm   \frac{dy^2}{Y(y)} \right)$& 
 \begin{minipage}{.3\textwidth}$-\left(\frac{\Im(h)}{\Im(h)^2 +\Re(h)^2}\right)^2dx^2 \\   +2\frac{\Re(h) \Im(h)}{   (\Im(h)^2 +\Re(h)^2)^2} dx dy  \\ +  \left(\frac{\Im(h)}{\Im(h)^2 +\Re(h)^2}\right)^2dy^2 $
 \end{minipage} &  \begin{minipage}{.3\textwidth}$  \frac{1+{x} Y'(y)}{Y(y)^4} \bigl(- 2Y(y) dxdy\\
    + (1+{x} Y'(y))dy^2\bigr)$\end{minipage}\\ \hline 
\end{tabular} 
where $h:=\Re(h) + i\cdot \Im(h)$ is  a holomorphic function of the variable $z:= x +i\cdot y $.   \end{Tha} 

\begin{rem} 
It it natural to consider the metrics from the Complex-Liouville case as the complexification of the metrics from the Liouville case: indeed, in the complex coordinates $z= x + i\cdot y$, $\bar z= x - i\cdot y$, the metrics  have the form

 \begin{equation*}  \left. \begin{array}{ccc}   ds^2_g & =  &  -\tfrac{1}{4}( \overline{h( z}) - h(z) )\left(d\bar z^2 -   dz^2\right),\\   
     ds^2_{\bar g} & =  &  -\tfrac{1}{4}\left(\frac{1}{\overline{h( z})} - \frac{1}{h(z)} \right)\left(\frac{d\bar z^2}{ \overline{h(z})} - \frac{ dz^2}{h(z)}\right).\end{array}\right. \end{equation*}
\end{rem} 

\begin{rem}  In the Jordan-block case, 
if $dY\ne 0$ (which is always the case  at almost every point,  if the restriction of  $g$ to any neighborhood   does   not admit a Killing vector field), after a  local coordinate change, the metrics  $g$ and   $\bar g$  have the form   \begin{eqnarray*}ds_g^2 & =&  \left( \tilde Y(y) +{x}{} \right)dxdy \\ 
   ds_{\bar g}^2   & =&  -\frac{2(\tilde Y(y)+x)}{y^3}dxdy  + \frac{(\tilde Y(y)+x)^2}{y^4}dy^2.
 \end{eqnarray*}
  \end{rem} 
  
We see that  the metric $g$ 
 from Complex-Liouville and  Jordan-block cases  always have signature $(+,-)$, and the metric $g$ 
 from  the Liouville case has  signature $(+,+)$ or $(-,-)$, if the sign ``$\pm$" is ``$+$".  
In this case, 
the formulas from Theorem A  are precisely the formulas obtained by Dini in \cite{Dini}.

We do not insist that we are the first to find these normal forms of projectively equivalent pseudo-Riemannian  metrics. According to \cite{Aminova},   a description of projectively equivalent metrics  was obtained by 
  P. Shirokov in \cite{shirokov}. Unfortunately, we were not able to find the reference \cite{shirokov} to check it.  The result of Theorem A could be     even more classical, see  Remark \ref{D}.

Given  two projectively  equivalent metrics,  it is easy to understand  what case they belong to. Indeed, the $(1, 1)$-tensor 
  $
  G^i_{j}:= \sum_{\alpha=1}^{2} \bar g_{j\alpha} g^{i\alpha}, 
  $   where $ g^{i\alpha}$ is inverse to $ g_{i \alpha}$,  
   has two different real eigenvalues in the Liouville  case, two complex-conjugated eigenvalues in the Complex-Liouvulle  case, and is (conjugate to) a Jordan-block in the Jordan-block case.

 There exists  an interesting and useful  connection of projectively equivalent metrics with integrable systems. 

Recall that 
 a function $F:{T^*}M^2\to \mathbb{R}$ is called \emph{an integral} of the geodesic flow of $g$, if $\{ H, F\}=0$, where $H:= \tfrac{1}{2} g^{ij} p_ip_j:T^*M^2\to \mathbb{R}$ is the kinetic energy corresponding to the metric, and $\{ \ , \ \}$ is the standard Poisson bracket on $T^*M^2$. Geometrically, this condition means 
that the function is constant on the orbits of the Hamiltonian system with the Hamiltonian $H$.  We say the  integral $F$ is {\it quadratic in  momenta,} if for  every local coordinate system $(x,y)$ on $M^2$  it has the form 
\begin{equation} \label{integral} F(x,y, p_x, p_y)= 
a(x,y)p_x^2+ b(x,y)p_xp_y+ c(x,y)p_y^2  
\end{equation} in the canonical coordinates $(x,y,p_x, p_y)$ on $T^*M^2$. 
Geometrically, the  formula \eqref{integral} 
means that the restriction of the integral to every cotangent space $T^*_pM^2\equiv \mathbb{R}^2$  is a homogeneous quadratic function. Of course, $H$ itself  is  an integral quadratic  in the momenta for $g$. We will say that the integral $F$ is {\it nontrivial}, if $F\ne  \textrm{const} \cdot H$ for all $\textrm{const} \in \mathbb{R}$. 

\begin{Thb} \label{main}  
 Suppose the metric $g$ on  ${M}^2$ admits a nontrivial integral quadratic in momenta. Then, in a neighbourhood of almost every point there exist coordinates  $(x,y)$ such that  the metric and the integral are as in the following table
\begin{center}\begin{tabular}{|c||c|c|c|}\hline & \textrm{Liouville case} & \textrm{Complex-Liouville case} & \textrm{Jordan-block case}\\ \hline \hline 
$g$ & $(X(x)-Y(y))(dx^2 \pm dy^2)$ &  $\Im(h)dxdy$ & $\left(1+{x} Y'(y)\right)dxdy $
\\ \hline  $F$&$\tfrac{X(x)p_y^2 \pm Y(y)p_{x}^2}{X(x)-Y(y)} $&  
 $p_x^2 - p_y^2   + 2\tfrac{\Re(h)}{\Im(h)}p_xp_y, $ & $p_x^2 -2  \frac{Y(y)}{1+ {x}{} Y'(y) }p_xp_y$  \\ \hline     
\end{tabular} \end{center}
$\textrm{where $h:=\Re(h) + i\cdot \Im(h)$ is a}$  \textrm{holomorphic function} $\textrm{of  the variable $z:= x+ iy$.}$
\end{Thb}

  Indeed, as it was shown in  \cite{MT,dim2},  and as it was essentially known to  Darboux~\cite[\S\S600--608]{Darboux}, if two metrics $g$ and $\bar g$ are projectively equivalent, then \begin{equation}\label{Integral}  I:TM^2\to \mathbb{R}, \
\
 I(\xi):=\bar g(\xi,\xi)\left(\frac{\det(g)}{\det(\bar g)}\right)^{2/3}
\end{equation}
 is an integral of the geodesic flow of  $g$. Moreover, it was shown in \cite{bryant}, see Section 2.4 there,  see also \cite{quantum},    the above statement  is proven to be  true\footnote{with a good will, one also can attribute this result to Darboux} in the  other direction: if the function \eqref{integral} is an integral for the geodesic flow of $g$,  then the metrics $g$ and $\bar g$ are projectively equivalent. Thus, Theorem A and Theorem B are equivalent. 
In this paper,  we will actually prove Theorem B obtaining Theorem A  as its consequence.

 \begin{rem} The corresponding natural Hamiltonian problem on the hyperbolic plane has  been recently treated in \cite{PR} following the approach used by Rosquist and Uggla \cite{ru:kt}. \end{rem} 

 \begin{rem} \label{D}  The formulas that will appear in the proof  are very close to that in \S 593 of \cite{Darboux}. Darboux  worked over complex numbers and therefore did  not care about whether the metrics are Riemannian or  pseudo-Riemannian. For example, Liouville and Complex-Liouville case are the same for him.  Moreover,   in \S 594, Darboux gets  
the  formulas that are very close to that of Jordan-block case, though he was interested  in the Riemannian case only, and, hence, treated this ``imaginary" case as not interesting.  \end{rem}

\section{Proof of Theorem B  (and, hence, of Theorem A)} 
If the metric $g$ has signature $(+, +)$ or $(-, -)$, Theorem A and, hence, Theorem B,  were obtained by Dini in \cite{Dini}. Below we assume that the metric $g$ has signature $(+, -)$.

\subsection{Admissible coordinate systems and 
Birkhoff-Kolokoltsov forms} \label{admissible}

Let $g$ be a pseudo-Riemannian metric on $M^2$  of signature $(+, -)$. 
Consider (and fix)  two linear independent  vector fields $V_1, V_2$ on $M^2$  such that 
\begin{itemize}  
 \item  $g(V_1, V_1) =g(V_2, V_2)=0$ and 
 \item $g(V_1, V_2)>0$.
 \end{itemize} 
Such vector fields always exist locally (and, since our result is local, this is sufficient for our proof). \weg{For possible further use, let us note that such vector fields always exist   on a finite (at most, 4-sheet-) cover of $M^2$.}

We will say that a local  coordinate system $(x,y)$
 is {\it admissible}, if  the vector fields  $\frac{\partial }{\partial x}$ and  $\frac{\partial }{\partial y} $ are proportional to $V_1, V_2$ with positive coefficient of proportionality:  $$V_1(x,y)=  \lambda_1(x,y)\frac{\partial }{\partial x} , \  \  \ V_2(x,y) =  \lambda_2(x,y) \frac{\partial }{\partial y}, \  \  \ \textrm{where $\lambda_i>0, \;\; i=1,2$}.$$ 
 Obviously, 
\begin{itemize} 
\item  admissible coordinates exist in a sufficiently small neighborhood  of every point, \item the metric $g$  in  admissible coordinates has the form 
 \begin{equation}\label{metric} 
 ds^2 =f(x,y)dxdy , \  \  \ \textrm{where $f>0$},  \end{equation}       
     \item two admissibe  coordinate systems  in 
      one neighbourhood  are connected by \begin{equation} \label{coordinatechange} \begin{pmatrix} x_{new}\\
 y_{new}\end{pmatrix} 
 = \begin{pmatrix} x_{new}(x_{old}) \\
 y_{new}(y_{old})\end{pmatrix} , \  \ \textrm{where  $\frac{dx_{ new}}{dx_{old}}>0$, $\frac{dy_{ new}}{dy_{old}}>0$}. \end{equation}   
\end{itemize}

\begin{Lemma}  \label{BK} 
Let $(x,y)$ be an admissible coordinate system for $g$. 
Let  $F$ given by  \eqref{integral} be  an  integral for  $g$. 
Then,    $B_1:= \frac{1}{\sqrt{|a(x,y)|}}dx$  ($B_2:= \frac{1}{\sqrt{|c(x,y)|}}dy$, respectively) is a    1-form, which is defined at points such that $a\ne 0$ ($c\ne 0$, respectively).  Moreover, the coefficient 
 $a$ ($c$, respectively) depends only on $x$ ($y$, respectively), which in particular imply that the forms $B_1$, $B_2$ are  closed.   \end{Lemma} 

\begin{rem}  The forms $B_1, B_2$ are not the direct analog of the ``Birkhoff" 2-form introduced by Kolokoltsov in \cite{Kol}. In a certain sense, they are  the real analog of the different branches of the  square root of the Birkhoff  form.   
\end{rem} 
{\bf Proof of the Lemma.} The first part of the statement, namely that the $\frac{1}{\sqrt{|a|}} dx$ ($\frac{1}{\sqrt{|c|}}dy$, respectively) transforms as a $1$-form under admissible coordinate changes is evident:  indeed, after the coordinate change  
\eqref{coordinatechange}, the  momenta transform as follows: 
 $p_{x_{old}}= p_{x_{new}}\frac{d{x_{new}}}{d{x_{old}}}$, $p_{x_{old}}= p_{x_{new}}\frac{d{x_{new}}}{d{x_{old}}}$. Then, the integral $F$ in the new coordinates has 
  the form $$ \underbrace{{a}\left(\frac{d{x_{new}}}{d{x_{old}}}\right)^2}_{a_{new} } {p_{x_{new}}^2}  + \underbrace{{b}\frac{d{x_{new}}}{d{x_{old}}}\frac{d{y_{new}}} {d{y_{old}}}}_{b_{new}}  {p_{x_{new}}} {p_{y_{new}}} + \underbrace{{c}\left(\frac{d{y_{new}}}{d_{y_{old}}} \right)^2}_{c_{new}} {p_{y_{new}}^2}.$$  Then, the formal expression $\frac{1}{\sqrt{|a|}}dx_{old}$
  $ \left(  \frac{1}{\sqrt{|c|}}dy_{old}, \textrm{ respectively}\right)$   transforms in  $$\frac{1}{\sqrt{|a|}} \frac{d{x_{old}}}{d{x_{new}}} dx_{new} \ \ \ \ \  \left(\textrm{  $ \frac{1}{\sqrt{|c|}}\frac{d{y_{old}}}{d{y_{new}}}dy_{new}$, respectively}\right),  $$ 
  which is precisely the transformation law of  1-forms.

Let us prove that the forms are closed. 
 If $g$ is given by \eqref{metric}, its Hamiltonian $H$ is given by $\frac{p_xp_y}{2f}$, and the condition $0=\{H, F\}$ reads \\
 $\begin{array}{ccl}
0&= &\left\{\frac{p_xp_y}{2f}, ap_x^2+ bp_xp_y+ cp_y^2\right\}  \\ &=& \frac{1}{f}\left(p_x^3(fa_y) + p_x^2 p_y (fa_x + fb_y + 2 f_x a + f_y b)+ p_yp_x^2 (fb_x + fc_y+ f_x b + 2 f_y)+ p_y^3 (c_xf)\right),\end{array}\\$ i.e., is equivalent to the following system of PDE: 
 \begin{equation}\label{sys} 
 \left\{\begin{array}{rcc} a_y&=&0 \\ fa_x + fb_y + 2 f_x a + f_y b&=&0\\  
 fb_x + fc_y+ f_x b + 2 f_yc &=&0\\ c_x&=&0 \end{array}
 \right.\end{equation} 
 
 Thus, $a=a(x)$, $c=c(y)$, which is equivalent to 
  $B_1:= \frac{1}{\sqrt{|a|}} dx$ and $B_2:= \frac{1}{\sqrt{|c|}}dy$ are closed forms (assuming $a\ne 0$ and $c\ne 0$). \qed

\begin{rem} \label{rem3} For further use let us formulate one more consequence of the equations \eqref{sys}: if $a\equiv c \equiv 0$ in a neighborhood of  a point, then $bf = \const$ implying     $F\equiv  \const \cdot H$ in the neighborhood.  
\end{rem}

Assume $a\ne 0$  ($c\ne 0$, respectively) 
 at a point $P_0$. For every point $P_1$ 
in a small neighbourhood $U$ of $P_0$ consider 
 \begin{equation} \label{normal}
 x_{new} :=\int\limits_{\begin{array}{c}
\gamma:[0,1]\to U \\ \gamma(0)=P_0 \\ \gamma(1)= P_1 \end{array}} 
 B_1, \ \ \left(\textrm{$ y_{new} :=\int\limits_{\begin{array}{c}
\gamma:[0,1]\to U \\ \gamma(0)=P_0, \\  \gamma(1)= P_1 \end{array}} 
B_2, $  respectively}\right).\end{equation}

Locally, in the admissible coordinates, 
 the functions $x_{new}$ and $y_{new}$ are given by  
 \begin{equation*}
  x_{new}(x{})=\int_{x_0}^{x_{}} \frac{1}{\sqrt{|a(t)|} }\, dt,  \ \ \  \   y_{new}(y)=\int_{y_0}^{y_{}} \frac{1}{\sqrt{|c(t)|} }\, dt. 
\end{equation*}

The new  coordinates $(x_{new}, y_{new})$ (or $(x_{new},y_{old})$ if $c_{old}\equiv 0$, 
or $\bigl((x_{old}, y_{new})$ if $a_{old}\equiv 0$ $\bigr)$    are admissible.   In these coordinates,   the forms $B_1$ and $B_2$     are given by  $\textrm{sign}(a_{old}) dx_{new}$,  $\textrm{sign}(c_{old})dy_{new}$  (we assume $\textrm{sign}(0)=0$).

\subsection{Proof of Theorem B}
We assume that $g$ of  signature $(+,-)$ on  $M^2$  
  admits  a nontrivial quadratic integral $F$ given by \eqref{integral}. Consider 
 the matrix $
 F^{ij}= \begin{pmatrix} a & \tfrac{b}{2} \\ \tfrac{b}{2} & c \end{pmatrix} $. It 
  can be viewed as  a $(2,0)$-tensor: if we change the coordinate system and rewrite the  function $F$ in the new coordinates, the matrix changes according to the tensor rule. Then, 
  \begin{equation*}
   \sum_{\alpha=1}^2 g_{j\alpha} F^{i\alpha}
  \end{equation*}  
  is a $(1,1)$-tensor. 
   In  a neighborhood $U$  of almost every point the Jordan normal form of this $(1,1)$-tensor is one of the following matrices: 
  
 {Case 1} $\begin{pmatrix} \lambda & 0\\ 0 & \mu  \end{pmatrix} $, \ \  \  \  \  \  {Case 2}   $\begin{pmatrix} \lambda+ i \mu  & 0\\ 0 & \lambda- i \mu   \end{pmatrix} $, \ \   \  \  \  
  {Case 3} $\begin{pmatrix} \lambda  & 1\\ 0 & \lambda   \end{pmatrix} $, \\
  where $\lambda,  \mu:U\to \mathbb{R}$.  
 Moreover,  in view of Remark~\ref{rem3}, there exists a neighborhood of almost every point such that $\lambda \ne \mu$ in Case 1 and $\mu \ne 0 $ in Case 2. 
  In the admissible coordinates, up to multiplication of $F $ by $-1$, and renaming $V_1\leftrightarrow V_2$,   
  Case 1 is equivalent to the condition 
    $a>0, \ c>0$,  Case 2 is equivalent to the condition   $a>0,$  $c<0$,  and Case 3 is equivalent to the   condition   $c\equiv 0$. 
  
  We now consider all three cases.
  
  \subsubsection{Case 1:  $a>0, \ c>0$. } 
 Consider the    coordinates \eqref{normal}. 
 In this coordinates,  $a=1$, $c=1$, and equations \eqref{sys} are: 
  
 \begin{equation*}
 \left\{\begin{array}{rcc}  (fb)_y+ 2 f_x  &=&0,\\  
 (fb)_x  + 2 f_y&=&0. \end{array}
 \right.\end{equation*} 
 This system can be solved. Indeed, it is equivalent to 
 \begin{equation*}
 \left\{\begin{array}{rcc}  (fb+ 2 f)_x  + (fb  + 2 f)_y&=&0,\\  
   (fb  -2 f)_x  -(fb- 2 f)_y&=&0, \end{array}
 \right.\end{equation*} 
 which,  after the change of cordinates $x_{new} = x+y$, $y_{new}= x-y$, has the form 
 \begin{equation*}
 \left\{\begin{array}{rcc}  (fb+ 2 f)_x&=&0,\\  
   (fb  -2 f)_y&=&0, \end{array}
 \right.\end{equation*} 
 implying $fb+ 2f = Y(y)$, $fb-2f=X(x)$. Thus, 
 $f= \tfrac{Y(y)-X(x)}{4}$, $b=2 \tfrac{X(x)+Y(y)}{Y(y)-X(x)}$.
 
 Finally, in the new coordinates, the metric and the integral have (up to a  possible multiplication by a constant) the form

 \begin{equation*}
 (X-Y)(dx^2 -dy^2)  
 \end{equation*}
 \and
 \begin{equation*}
 2\left(p_x^2 -  \tfrac{X(x)+Y(y)}{X(x)-Y(y)}(p_x^2-p_y^2) + p_y^2\right)= 4\frac{p_y^2 X(x) - p_x^2Y(y)}{X(x)-Y(y)}.  
 \end{equation*}
 Theorem B is proved under the assumptions of Case 1.

  \subsubsection{ Case 2:   
  $a>0$, $c<0$. }

 Consider the coordinates \eqref{normal}. 
 In this coordinates,  $a=1$, $c=-1$, and the equations \eqref{sys} are:
 
 \begin{equation*}
 \left\{\begin{array}{rcc}  (fb)_y+ 2 f_x  &=&0,\\  
 (fb)_x  -2 f_y&=&0. \end{array}
 \right.\end{equation*} 
 We see that these conditions are the Cauchy-Riemann  conditions for the complex-valued 
 function $fb+ 2i \cdot f$. Thus, for an appropriate holomorphic 
 function $h= h(x+ i\cdot y)$,  we have $fb=\Re(h)$, $2f =\Im(h)$. 
 Finally, in a certain coordinate system the metric and the integral are 
 (up to multiplication by constants): 
  \begin{equation*}
 2\Im(h)dxdy \ \ \ \textrm{and} \ \ \ p_x^2 - p_y^2   + 2\tfrac{\Re(h)}{\Im(h)}p_xp_y.   
 \end{equation*}
 Theorem B is proved under the assumptions of Case 2.

 \subsubsection{Case 3: $a>0$, $c=0$.  } 
 
 Consider  admissible 
 coordinates $x, y$, such that  $x$ is 
 the coordinate from  \eqref{normal}. 
 In these coordinates,   $a=1$, $c=0$, and the equations \eqref{sys} are: 
 \begin{equation*}
 \left\{\begin{array}{rcc}  (fb)_y+ 2 f_x  &=&0\\  
 (fb)_x  &=&0 \end{array}.
 \right.\end{equation*} 
 This system can be solved. Indeed, the second equation implies $fb= -Y(y)$. Substituting this in the first equation  we obtain 
 $Y'= 2f_x$ implying $$f= \frac{x}{2} Y'(y)+ \widehat{ Y}(y) \textrm{ \ \ and}  \ \ \ b= - \frac{Y(y)}{\frac{x}{2} Y'(y)+ \widehat{ Y}(y)}.$$ Finally, the metric and the integral are 
  \begin{equation}
  \label{answer:case3} 
 \left( \widehat{ Y}(y)+\frac{x}{2} Y'(y)\right)dxdy \ \ \ \textrm{and} \ \ \ p_x^2 - \frac{Y(y)}{\widehat{ Y}(y)+\frac{x}{2} Y'(y) }p_xp_y
 \end{equation}

 Moreover, by the change $y_{new}=\beta(y_{old})$ the metric and the integral  \eqref{answer:case3} will be  transformed to:
  \begin{equation*}
 \left( \widehat{ Y}(y)\beta' +\frac{x}{2} Y'(y)\right)dxdy \ \ \ \textrm{and} \ \ \ p_x^2 +  \frac{Y(y)}{\widehat{ Y}(y)\beta'+\frac{x}{2} Y'(y) }p_xp_y
 \end{equation*}
 Thus, by putting $\beta(y) = \int_{y_0}^y \frac{1}{\widehat{ Y}(t)} dt$,  we can make the metric and the integral to be 
 $$  \left( 1+\frac{x}{2} Y'(y)\right)dxdy \ \ \ \textrm{and} \ \ \ p_x^2 -   \frac{Y(y)}{1+\frac{x}{2} Y'(y) }p_xp_y.
 $$
 Moreover, after the coordinate change $x_{new} =  \tfrac{x_{old}}{2}$ and dividing/multiplication of the metric/integral by ${2}$, the metric and the integral have  the form  from Theorem B
  \label{case3} 
  \begin{equation*}
    \left( 1+ {x} Y'(y)\right)dxdy \ \ \ \textrm{and} \ \ \ p_x^2 -2  \frac{Y(y)}{1+ {x}{} Y'(y) }p_xp_y
 \end{equation*}

 Theorem B is proved.
  


\begin{thebibliography}{99}

\bibitem{Aminova}
A. V. Aminova,
\emph{Projective transformations of pseudo--Riemannian manifolds. Geometry, 9.}
 J. Math. Sci. (N. Y.) \textbf{113} (2003), no. 3, 367--470.
 
\bibitem{Beltrami}
E. Beltrami,
\emph{Resoluzione del problema: riportari i punti di una
superficie sopra un piano in modo che le linee geodetische vengano
rappresentante da linee rette},
Ann. Mat., \textbf{1} (1865), no. 7, 185--204.


\weg{\bibitem{3} { G. D.  Birkhoff, }
\emph{Dynamical Systems, } A.M.S.  Colloq. Publ. {\bf 9}, Amer. Math. Soc., New York, 1927.} 



\bibitem{bryant}   R. L. Bryant,
G. Manno,
V.  S. Matveev, \emph{A solution of a problem of Sophus Lie: Normal forms of 2-dim metrics admitting two projective vector fields}, accepted to Math. Ann, 
   	arXiv:0705.3592 .   
   	
\bibitem{Darboux}
G. Darboux,
\emph{Le\c{c}ons sur la th\'eorie g\'en\'erale des surfaces},
Vol. III, Chelsea Publishing, 1896.


\bibitem{Dini} U. Dini, \emph{ Sopra un problema che si presenta nella teoria
generale delle rappresentazioni geografiche di una superficie su un'altra},
Ann. Mat., ser.2, {\bf 3}(1869), 269--293.


\bibitem{Kol} {V.~N. Kolokoltsov, } \emph{
Geodesic flows on two-dimensional manifolds with an
additional first integral that is polynomial with respect to
velocities, }
Math. USSR-Izv. {\bf 21}(1983), no.~2, 291--306. 

\bibitem{MT}
V. S. Matveev, P. J. Topalov,
\emph{Trajectory equivalence and corresponding integrals},
 Regular and Chaotic Dynamics,
\textbf{3} (1998), no. 2, 30--45.

\bibitem{dim2}
V. S. Matveev, P. J. Topalov,
\emph{Geodesic equivalence of metrics on surfaces, and their
integrability,} Dokl. Math. \textbf{60} (1999), no.1, 112-114.

\bibitem{Quantum} V.S.  Matveev, \emph{ Quantum integrability of the Beltrami-Laplace operator for geodesically equivalent metrics.}    Dokl. Akad. Nauk  {\bf 371}(2000),  no. 3, 307--310. 
  
  \bibitem{quantum}  V. S. Matveev, P. J. Topalov,
 {\it Quantum integrability for the   Beltrami-Laplace
 operator as geodesic  equivalence,}
  Math. Z. {\bf 238}(2001), 833--866.

\bibitem{PR}
G. Pucacco, K. Rosquist, \emph{(1+1)-dimensional separation of variables}, J. Math. Phys. (2007), {\bf 48}, 112903--112925.

\bibitem{ru:kt}
K. Rosquist, C. Uggla, \emph{Killing tensors in two-dimensional space-times with applications to cosmology}, J. Math. Phys.  (1991), {\bf 32},  3412--3422.

\bibitem{shirokov} P. A. Shirokov, \emph{Selected Works on Geometry}, Kazan Univ., Kazan (1966), 383--389.
\weg{\bibitem{Whi}
{ Whittaker, E. T., }
\emph{ A Treatise on the Analytical Dynamics of Particles and Rigid
  Bodies, }
Cambridge University Press, Cambridge, 1937.}
  
\end{thebibliography}
 \end{document}